\documentclass[12pt]{amsart}
\usepackage{amssymb,amsfonts,amsthm}

\newcommand{\bi}{\bibitem}
\newcommand{\nb}{\newblock}

\newcommand{\zz}{{\mathbb{Z}}}
\newcommand{\be}[1]{\begin{equation}\label{#1}}
\newcommand{\ee}{\end{equation}}

\newtheorem{thm}{\quad Theorem}
\newtheorem{lm}{\quad Lemma}

\newtheorem*{question}{Question}

\let\cal\mathcal
\let\om\omega

\begin{document}

\title
{Growth rates of amenable groups}
\author{G.~N.~Arzhantseva}
\address{Section de Math\'ematiques\\
    Universit\'e de Gen\`eve\\ CP 240, 1211 Gen\`eve 24\\
    Switzerland}
\email{Goulnara.Arjantseva@math.unige.ch}
\author{V.~S.~Guba}
\address{Department of mathematics\\
    Vologda State University\\ 6 S. Orlov St.\\ Vologda\\ 160600\\
    Russia}
\email{guba@uni-vologda.ac.ru}
\author{L.~Guyot}
\address{Section de Math\'ematiques\\
    Universit\'e de Gen\`eve\\ CP 240, 1211 Gen\`eve 24\\
    Switzerland}
\email{Luc.Guyot@math.unige.ch}

\thanks{The work has been supported by the Swiss National Science
 Foundation, No.~PP002-68627.}

\date{}

\begin{abstract}
Let $F_m$ be a free group with $m$ generators and let $R$ be its normal
subgroup such that $F_m/R$ projects onto $\zz$. We give a lower bound for the
growth rate of the group $F_m/R'$ (where $R'$ is the derived subgroup of $R$)
in terms of the length $\rho=\rho(R)$ of the shortest nontrivial relation in
$R$. It follows that the growth rate of $F_m/R'$ approaches $2m-1$ as $\rho$
approaches infinity. This implies that the growth rate of an $m$-generated
amenable group can be arbitrarily close to the maximum value $2m-1$. This
answers an open question by P. de la Harpe. In fact we prove that such groups
can be found already in the class of abelian-by-nilpotent groups as well as in
the class of finite extensions of metabelian groups.
\end{abstract}

\maketitle

\section{Introduction}

Let $G$ be a finitely generated group and $A$ a fixed finite set of
generators for ~$G$. By $\ell(g)$ we denote the {\em word length\/}
of an element $g\in G$ in the generators ~$A$, i.e.\ the length of
a shortest word in the alphabet ~$A^{\pm 1}$ representing ~$g$. Let
$B(n)$ denote the ball $\{g\in G\mid\ell(g)\le n\}$ of radius ~$n$
in ~$G$ with respect to ~$A$. The {\it growth rate\/} of the pair
$(G,A)$ is the limit
$$
\om(G,A)=\lim_{n\to\infty}\sqrt[n]{|B(n)|}.
$$
(Here $|X|$ denotes the number of elements of a finite set ~$X$.)
This limit exists due to the submultiplicativity property of the
function ~$|B(n)|$,
see for example \cite[VI.C, Proposition 56]{Harpe}. Clearly,
$\om(G,A) \ge 1$. A finitely generated group ~$G$ is said to be of
{\em exponential growth\/} if $\om(G,A)>1$ for some (which in fact
implies for any) finite generating set ~$A$. Groups with $\om(G,A)=1$
are groups of {\em subexponential\/} growth.

Let $|A|=m$. It is known that $\om(G,A)=2m-1$ if and only if
$G$ is freely generated by $A$~\cite[Section V]{GrH}.
In this case $G$ is non-amenable whenever $m>1$.

A finitely generated group which is nonamenable is necessarily of
exponential growth~\cite{Adel}. The following interesting question
is due to P.~de la Harpe.

\begin{question}\cite[VI.C~62]{Harpe}
For an integer $m\ge2$, does there exist a constant $c_m>1$, with
$c_m<2m-1$, such that $G$ is not amenable provided  $\om(G,A)\ge c_m$?
\end{question}

We show that the answer to this question is negative. Thus, given $m\ge2$,
there exists an amenable group on $m$ generators with the growth rate
as close to $2m-1$ as one likes.

It is worth noticing that for every $m\ge2$ there exists a sequence of
non-amenable groups (even containing  non-abelian free subgroups) whose
growth rates approach 1 (see~\cite{GrHa}).
\vspace{0.5ex}

For a group $H$, we denote by $H'$ its derived subgroup, that is,
$[H,H]$.
\vspace{0.5ex}

The authors thank A.\,Yu.\,Ol'shanskii for helpful comments.

\section{Results}

Let $F_m$ be a free group of rank $m$ with free basis $A$. Suppose
that $R$ is a normal subgroup of $F_m$. Assume that there is a
homomorphism $\phi$ from $F_m$ {\bf onto} the (additive) infinite cyclic
group such that $R$ is contained in its kernel (that is, $F_m/R$ has $\zz$
as a homomorphic image). By $a$ we denote a letter from $A^{\pm1}$ such that
$$
\phi(a)=\max\{\,\phi(x)\mid x\in A^{\pm1}\,\}.
$$
Clearly, $\phi(a)\ge1$.
\vspace{1ex}

Throughout the paper, we fix a homomorphism $\phi$ from $F_m$ {\bf onto}
$\zz$, the letter $a$ described above and the value $C=\phi(a)$. By $R$
we will usually denote a normal subgroup in $F_m$ that is contained in
the kernel of $\phi$.
\vspace{1ex}

A word $w$ over $A^{\pm1}$ is called {\em good\/} whenever it satisfies the
following conditions:

\begin{enumerate}
\item
$w$ is freely irreducible,
\item
the first letter of $w$ is $a$,
\item
the last letter of $w$ is not $a^{-1}$,
\item
$\phi(w)>0$.
\end{enumerate}

Let $D_k$ be the set of all good words of length $k$ and let $d_k$ be
the number of them.

\begin{lm}
\label{estim}
The number of good words of length $k\ge4$ satisfies the following inequality:
\be{goods}
d_k\ge4m(m-1)^2(2m-1)^{k-4}.
\ee
In particular, $\lim\limits_{k\to\infty}d_k^{1/k}=2m-1$.
\end{lm}

\proof
Let $\Omega$ be the set of all freely irreducible words $v$ of length $k-1$
satisfying $\phi(v)\ge0$. The number of all freely irreducible words of length
$k-1$ equals $2m(2m-1)^{k-2}$. At least half of them has a nonnegative image
under $\phi$. So $|\Omega|\ge m(2m-1)^{k-2}$.

Let $\Omega_1$ be the subset of $\Omega$ that consists of all words whose
initial letter is different from $a^{-1}$. We show that
$|\Omega_1|\ge((2m-2)/(2m-1))|\Omega|$. It is sufficient to prove that
$|\Omega_1\cap A^{\pm1}u|\ge((2m-2)/(2m-1))|\Omega\cap A^{\pm1}u|$ for any
word $u$ of length $k-2$. Suppose that $a^{-1}u$ belongs to $\Omega$.
For every letter $b$ one has $\phi(b)\ge\phi(a^{-1})$. Therefore,
$bu\in\Omega_1$ for every letter $b\ne a^{-1}$ provided $bu$ is irreducible.
There are exactly $2m-2$ ways to choose a letter $b$ with the above properties.
Hence $|\Omega_1\cap A^{\pm1}u|$ and $|\Omega\cap A^{\pm1}u|$ have $2m-2$
and $2m-1$ elements, respectively. If $a^{-1}u\notin\Omega$, then both sets
coincide.

Now let $\Omega_2$ denote the subset of $\Omega_1$ that consists of all
words whose terminal letter is different from $a^{-1}$. Analogous argument
implies that $|\Omega_2|\ge((2m-2)/(2m-1))|\Omega_1|$. It is obvious that
$av$ is good provided $v\in\Omega_2$. Therefore, the number of good words
is at least
$$
|\Omega_2|\ge\frac{2m-2}{2m-1}|\Omega_1|\ge\left(\frac{2m-2}{2m-1}\right)^2
|\Omega|\ge4m(m-1)^2(2m-1)^{k-4}.
$$
\endproof

To every word $w$ in $A^{\pm1}$ one can uniquely assign the path $p(w)$ in the
Cayley graph ${\cal C}={\cal C}(F/R,A)$ of the group $F/R$ with $A$ the
generating set. This is the path that has label $w$ and starts at the identity.
We say that a path $p$ is {\em self-avoiding\/} if it never visits the same
vertex more than once.

Let $\rho=\rho(R)$ be the length of the shortest nontrivial element in a
normal subgroup $R\le F_m$.

\begin{lm}
\label{saw}
Let $R$ be a normal subgroup in $F_m$ that is contained in the kernel of a
homomorphism $\phi$ from $F_m$ onto $\zz$. Suppose that $k\ge2$ is chosen in
such a way that the following inequality holds:
\be{rhok}
\rho(R)>Ck(2k-3)+2k-2.
\ee
Then any path in the Cayley graph ${\cal C}$ of $F_m/R$ labelled by a word of
the form $g_1g_2\cdots g_t$, where $t\ge1$, $g_s\in D_k$ for all $1\le s\le t$,
is self-avoiding.
\end{lm}

\proof
If $p$ is not self-avoiding, then let us consider its minimal subpath $q$
between two equal vertices. Clearly, $|q|\ge\rho\ge k$. Therefore, $q$ can be
represented as $q=g'g_i\cdots g_jg''$, where $g_i$, \dots, $g_j$ are in $D_k$,
the word $g'$ is a proper suffix of some word in $D_k$, the word $g''$ is a
proper prefix of some word in $D_k$. We have $|g'|,|g''|\le k-1$ so
$|g_i\ldots g_j|>Ck(2k-3)$. This implies that $j-i+1$ (the number of sections
that are completely contained in $q$) is at least $C(2k-3)+1$. Obviously,
$\phi(g')\ge-C(k-1)$ and $\phi(g'')\ge-C(k-2)$ (we recall that $g''$ starts
with $a$ if it is nonempty). On the other hand, $\phi(g_s)\ge1$ for all $s$.
Hence $\phi(g_i\cdots g_j)\ge j-i+1\ge C(2k-3)+1$ and so
$\phi(g'g_i\cdots g_jg'')\ge1$, which is obviously impossible because for every
$r\in R$ one has $\phi(r)=0$.
\endproof

\begin{thm}
\label{main}
Suppose that $R$ is a normal subgroup in the free group $F_m$ that is
contained in the kernel of a homomorphism $\phi$ from $F_m$ onto $\zz$.
Let $C$ be the maximum value of $\phi$ on the generators or their inverses.
Let $\rho=\rho(R)$ be the length of the shortest cyclically irreducible
nonempty word in $R$. If a number $k\ge4$ satisfies the inequality
\be{rhock}
\rho\ge Ck(2k-3)+2k-1,
\ee
then the growth rate of the group $F_m/R'$ w.r.t. the natural generators
is at least
$$
(2m-1)\cdot\left(\frac{4m(m-1)^2}{(2m-1)^4}\right)^{1/k}.
$$
\end{thm}

\proof
We use the following known fact \cite[Lemma 1]{DLS}. A word $w$ belongs to $R'$
if and only if, for any edge $e$, the path labelled by $w$ in the Cayley graph
of the group $F_m/R$ has the same number of occurrences of $e$ and $e^{-1}$. Hence,
if we have a number of different self-avoiding paths of length $n$ in the Cayley
graph of $F_m/R$, then they represent different elements of the group $F_m/R'$.
Moreover, all the corresponding paths in the Cayley graph of $F_m/R'$ are geodesic
so these elements have length $n$ in the group $F_m/R'$.

Suppose that the conditions of the theorem hold. For every $n$, one can
consider the set of all words of the form $g_1g_2\ldots g_n$, where all
the $g_i$'s belong to $D_k$. By Lemma \ref{saw} all these elements give
us different self-avoiding paths in the Cayley graph of $F_m/R$. Hence
for any $n$ we have at least $d_k^n$ different elements in the group $F_m/R'$
that have length $kn$. Therefore, the growth rate of $F_m/R'$ is at least
$d_k^{1/k}$. It remains to apply Lemma \ref{estim}.
\endproof

One can summarize the statement of Theorem \ref{main} as follows: if
all relations of $F_m/R$ are long enough, then the growth rate of
the group $F_m/R'$ is big enough. Notice that we cannot avoid the
assumption that $F_m/R$ projects onto $\zz$. Indeed, for any number
$\rho$, there exists a finite index normal subgroup in $F_m$ such
that all the nontrivial elements in this subgroup are longer than $\rho$.
If $R$ was such a subgroup, then $F/R'$ would be a finite extension of an
abelian group and its growth rate would be equal to $1$.

\begin{thm}
\label{chain}
Let $F_m$ be a free group of rank $m$ with free basis $A$ and
let $\phi$ be a homomorphism from $F_m$ onto $\zz$. Suppose that
$$
\ker\phi\ge R_1\ge R_2\ge\cdots\ge R_n\ge\cdots
$$
is a sequence of normal subgroups in $F_m$. If the intersection of all the
$R_n$'s is trivial, then the growth rates of the groups $F_m/R_n'$ approach
$2m-1$ as $n$ approaches infinity, that is,
$$
\lim\limits_{n\to\infty}\,\om(F_m/R_n',A)=2m-1.
$$
\end{thm}

\proof
Since the subgroups $R_n$ have trivial intersection, the lengths of their
shortest nontrivial relations approach infinity, that is,
$\rho(R_n)\to\infty$ as $n\to\infty$. Let $k(n)=\left[\sqrt{\rho(R_n)/2C}\right]$,
where $C$ is defined in terms of $\phi$ as above). Obviously, inequality
(\ref{rhock}) holds and $k(n)\to\infty$. Now Theorem \ref{main} implies that
the growth rates of the groups $F_m/R_n'$ approach $2m-1$.
\endproof

Now we show that for every $m$ there exists an amenable group with $m$
generators whose growth rate is arbitrarily close to $2m-1$.

\begin{thm}
\label{amen}
For every $m\ge1$ and for every $\varepsilon>0$, there exists an $m$-generated
amenable group $G$, which is an extension of an abelian group by a nilpotent
group such that the growth rate of $G$ is at least $2m-1-\varepsilon$.
\end{thm}

\proof
It suffices to take the lower central series in the statement of Theorem \ref{chain}
(that is, $R_1=F_m'$, $R_{i+1}=[R_i,F_m]$ for all $i\ge1$). The subgroups $R_n$ ($n\ge1$)
have trivial intersection and they are contained in $F_m'$. So all of them are
contained in the kernel of a homomorphism $\phi$ from $F_m$ onto $\zz$ (one of the
free generators of $F_m$ is sent to 1, the others are sent to 0). The groups
$G_n=F_m/R_n'$ are extensions of (free) abelian groups $R_n/R_n'$ by (free) nilpotent
groups $F_m/R_n$ so all the groups $G_n$ are amenable. The growth rates of them approach
$2m-1$.
\endproof

Notice that one can take the sequence $R_n=F_m^{(n)}$ of the $n$th derived subgroups
as well (that is, $R_1=F_m'$, $R_{i+1}=R_i'$ for all $i\ge1$). It is not hard to
show that $\rho(R_n)$ grows exponentially. The groups $F_m/R_n'=F_m/R_{n+1}$ are free
soluble. Their growth rates approach $2m-1$ very quickly. For instance, the growth rate
of the free soluble group of degree 15 with 2 generators is greater than $2.999$.

One more application of Theorem \ref{amen} can be obtained as follows. The
set of finite index subgroups of $F_m$ is countable so one can enumerate them as
$N_1$, $N_2$, \dots, $N_i$, \dots\,. Let $M_i=N_1\cap N_2\cap\cdots\cap N_i$
and let $R_i=M_i'$ for all $i\ge1$. Obviously, the subgroups $M_i$ (and thus $R_i$)
have trivial intersection. Indeed, $F_m$ is residually finite and so the subgroups
$N_i$ intersect trivially. As above, all the $R_i$'s are contained in $F_m'$ so
they are contained in the kernel of a homomorphism $\phi$ from $F_m$ onto $\zz$.
Hence the growth rates of the groups $F_m/R_i'=F_m/M_i''$ approach $2m-1$. These
groups are extensions of $M_i/M_i''$ by $F_m/M_i$, that is, they are finite
extensions of (free) metabelian groups.

Therefore, in each of the two classes of groups: 1) extensions of abelian
groups by nilpotent groups, 2) finite extensions of metabelian groups, there
exist $m$-generated groups with growth rates approaching $2m-1$.
\vspace{0.5ex}

{\bf Remark.}\ A.\,Yu.\,Ol'shan\-skii suggested the following improvement.
Let $p$ be a prime. Since $F_m$ is residually a finite $p$-group, one can get
a chain $M_1\ge M_2\ge\cdots$ of normal subgroups with trivial intersection, where
$F_m/M_i$ are finite $p$-groups. Now let $R_i=\ker\phi\cap M_i$. The group
$F_m/R_i$ is a subdirect product of $\zz$ and a finite $p$-group. In particular,
it is nilpotent. Besides, it is an extension of $\zz$ by a finite $p$-group
and an extension of a finite $p$-group by $\zz$, as well. So $F_m/R_i'$ will be
abelian-by-nilpotent and metabelian-by-finite at the same time. (In fact,
the metabelian part is an extension of an abelian group by $\zz$.) Also
one can view $F_m/R_i'$ as an extension of a virtually abelian group by $\zz$.


\begin{thebibliography}{10}
\bi{Adel} G.~M.~Adel'son-Vel'skii and Yu.~A.~Sreider,
{\em The Banach mean on groups}, Uspehi Mat. Nauk (N.S.) 12, 1957 no.~{\bf
6}(78), 131--136.


\bi{DLS} C.~Droms, J.~Lewin, H.~Servatius.
{\em The length of elements in free solvable groups},
Proc. Amer. Math. Soc. 119, {\bf 1} (1993), 27--33.

\bibitem{GrH}
R.~Grigorchuk and P.~de la Harpe,
{\em On problems related to growth, entropy, and spectrum in group theory},
J. Dynam. Control Systems {\bf 3} (1) (1997), 51--89.

\bibitem{GrHa} R.~Grigorchuk and P.~de la Harpe,
{\em Limit behaviour of exponential growth rates for finitely generated groups},
\nb Essays on geometry and related topics, Vol. 1, 2, 351--370,
Monogr. Enseign. Math., {\bf 38}, Enseignement Math., Geneva, 2001.

\bibitem{Harpe}
P.~de la Harpe, {\em Topics in geometric group theory}, Chicago Lectures in
Mathematics. University of Chicago Press, Chicago, IL, 2000.

\end{thebibliography}
\end{document}